# NECO - A scalable algorithm for NEtwork COntrol

Sean P. Cornelius[1], Adilson E. Motter[1,2]

[1] *Department of Physics and Astronomy, Northwestern University, Evanston, IL 60208, USA.*
[2] *Northwestern Institute on Complex Systems, Northwestern University, Evanston, IL 60208, USA.*

We present an algorithm for the control of complex networks and other nonlinear, high-dimensional dynamical systems. The computational approach is based on the recently-introduced concept of *compensatory perturbations*—intentional alterations to the state of a complex system that can drive it to a desired target state even when there are constraints on the perturbations that forbid reaching the target state directly. Included here is ready-to-use software that can be applied to identify eligible control interventions in a general system described by coupled ordinary differential equations, whose specific form can be specified by the user. The algorithm is highly scalable, with the computational cost scaling as the number of dynamical variables to the power 2.5.

**Keywords:** disordered systems, network control, cascades, transient stability, nonlinear dynamics

**Introduction:**

Whether natural or man-made, real complex networks are often characterized by the following properties:

- the dynamics is nonlinear;
- the system has multiple stable states (or attractors);
- the system is described by a large number of dynamical variables;
- there are constraints on the physically implementable control interventions;
- the network exhibits decentralized response to perturbations.

The final point, in particular, means that the system generally responds suboptimally to damage caused by external perturbations. Perturbations can propagate through the network, influencing components beyond those initially affected. In the recent associated publication by Cornelius et al. [1], it is shown that the same phenomenon can actually be used to design compensatory perturbations that can control a network. These interventions take advantage of the nonlinearity inherent to real networks, driving the system to a desired target state indirectly by bringing it to the so-called *basin of attraction* of that state. This principle allows one to control a broad class of networks; indeed, situations where the target is not directly reachable are expected to be the rule rather than the exception.

In this protocol we present NECO, an algorithm and software for NEtwork COntrol that allows systematic identification of such control interventions.

NECO implements a scalable method—originally introduced in ref. [1]—for identifying compensatory perturbations in general networks of dynamical units under arbitrary constraints. The algorithm itself is depicted in flowchart form in Fig. 1, and is based on iteratively updating the initial state of the network with small, incremental perturbations, each of which is forecasted to move the system's trajectory closer to the target state while at the same time complying with the given constraints. This procedure is terminated once an eligible state that evolves to the target is found.

**Equipment:**

The source code provided with this protocol implements the control algorithm introduced in ref. [1] as a Python module called *neco*. The code is documented below and requires only the following freely-available software/libraries:

-Python (http://www.python.org/)
-SciPy (http://www.scipy.org/)
-NumPy (http://www.numpy.org/)

**Procedure:**

The complex networks we envision controlling are the common case of nonlinear, dissipative dynamical systems, where the time evolution of the component dynamical variables is governed by ordinary differential equations (ODEs). The roles of the individual nodes and the structure of the network as a whole are assumed to be reflected in these ODEs. For example, if there is a link in the network from node *i* to node *j*, one expects a term in the *j*-th differential equation involving the *i*-th state variable. Thus for the operational purpose of "control", we can completely characterize a given complex network by:
- a (vector) function "f" that gives the "right-hand side" of the first-order ODEs defining the dynamics;
- two (vector) functions "g" and "h" that define, respectively, the inequality and equality constraints on the eligible interventions.

Note that from this point onward, we are no longer talking about NECO in the abstract, but rather in regard to its software implementation in Python. Thus f, g, and h here denote Python functions that are to be implemented by the user and input to *neco*. These should be regarded as distinct from the corresponding mathematical abstractions in Fig. 1, which are similarly labeled as F, g, and h (in Roman font) to facilitate comparison. In particular, the Python functions will generally take different arguments, and in a different order, than their abstract counterparts. For example, f—the Python implementation of the derivatives function F—can accept a number of optional parameters that are needed to compute the right-hand side of the ODEs, parameters that were omitted for notational clarity in Fig. 1.

Also, note that within the source code and from this point onward in this protocol, we use the letter 'y' instead of 'x' to denote the state of a dynamical system, to follow the library SciPy's notational conventions for systems of ODEs.

**Defining the ODE system.** To apply this software to find compensatory perturbations in a given dynamical system, one needs only to supply a Python function of the form

f(y, t, param1, param2, …),

which, given a system state y and optional parameters param1, param2, etc., evaluates the right-hand side of the ordinary differential equations that govern the time evolution of the system. We anticipate that the user will frequently want to use this function for other purposes, e.g., integration of system orbits for the purposes of visualization. Thus, we have intentionally required the call signature of f to take the form above so that it can be used with SciPy's tools for ODEs (in particular their numerical integrator, odeint) without modification. Note, however, that *neco* assumes that the ODE system described by f be autonomous—that is, it has no explicit time dependence. Hence, even though the time variable t appears as a required second argument, it should not actually be used within the calculation of the return value of f. The condition of being autonomous does not sacrifice any generality, since every non-autonomous system of ODEs (those with an explicit dependence on t) can be transformed into an equivalent autonomous system by augmenting the state space with an additional "dummy" variable representing t, whose time derivative is then the constant 1.

Note that we impose no specific conditions on the nature of the dynamics defined by f; notably, it can contain arbitrary nonlinearities and be arbitrarily high-dimensional. The only (mild) requirement is that f be differentiable with respect to the system state y, so that a Jacobian matrix can be calculated.

**Defining constraints on the eligible interventions.** Given the network state y0 at a particular time, one can specify the subset of network states y that are reachable through eligible perturbations from y0 using vector constraint functions of the form g(y) and h(y) for inequality and equality constraints, respectively. Each of these functions operates by taking a state y of the network and returning a NumPy array with lengths equal to the number of inequality or equality constraints, respectively. States that can be reached by eligible control perturbations are then those that satisfy g(y) <= 0 and h(y) == 0, where the inequality and equality are taken to apply individually to every element of the returned arrays.

A common occurrence of inequality/equality constraints is one in which some dynamical variables (components of y) are allowed to be perturbed over a certain range and/or cannot be perturbed at all. This special case of component-wise constraints can always be expressed using g and h above, but, for convenience, we have provided the option to specify vectors of lower (lb) and upper (ub) bounds on eligible perturbed states y, which must then satisfy

lb[i] <= y[i] <= ub[i],

for every component *i*. This includes in particular the case in which a state variable *i* cannot be perturbed at all, *i.e.,* lb[i] = ub[i] = y0[i].

Each of lb, ub, g, and h is an optional keyword argument for the function *neco*. If either vector lb or vector ub is not supplied, the coordinates of the eligible states are assumed to be unbounded from below or above, respectively, before the imposition of any constraints specified by g(y) and/or h(y). Individual components can be made unbounded with respect to lb and ub by filling the appropriate positions in these vectors with the NumPy constants -inf and +inf, respectively.

**Parameters for the control procedure.** The enclosed implementation of *neco* accepts a handful of parameters that govern its internal operation. For the motivation behind each of these parameters, as well as guidance in how to choose their values based on the time and length scales of a particular dynamical system, we refer the reader to ref. [1]. Each of these parameters is implemented as a keyword argument in the *neco* Python function, which we summarize here:

eps0:   Minimum size of the incremental perturbation at every iteration (eps0 > 0), to be applied at the initial time to ensure that the algorithm makes a non-zero step every time.
eps1:   Maximum size of the incremental perturbation at every iteration (eps1 > eps0), to ensure that the forecasted perturbation using the variational matrix is valid (at both the initial time and at the forecasted time of closest approach to the target).
it_max: Maximum number of iterations of the algorithm. If this number is exceeded before finding a compensatory perturbation, the function will return 1, indicating failure.
t_max:  Time window over which to search for the closest approach to the target at every iteration.
dt:     Integration time step.
t_test: Time window over which to test convergence of the system's orbit to the target state (generally t_test >> t_max)
tol:    Numerical tolerance for convergence to the target state.

Thus, if a state y evolves to within a ball of radius tol around the target state, and within t_test time units, it is considered to be in the basin of attraction of the target and the algorithm declares success.

**Other optional parameters.**
jac(y, t, param1, param2, …): Function that computes the Jacobian matrix J of f as a 2D NumPy array, where J[i, j] = df[i]/dx[j]. If not supplied, the Jacobian will be calculated from f numerically.
dist(y1, y2):   Function that defines the distance between two states, y1 and y2, in the system under consideration. By default, the Euclidean norm is used.
n_test:         Every n_test iterations, test whether the current initial state attracts to the target. By default, n_test = 1.

full_output: If full_output is False, the method simply returns the final perturbed state. If full_output is True, the method returns a dictionary containing more detailed information (see "Return values" below). This option is False by default.

args: A tuple of additional parameters (other than y and t) that will be passed to the derivatives function f and the Jacobian function jac (if supplied).

**Return values.**
(if full_output is False)
y0_prime: The final perturbed state.

(if full_output is True)
info: Dictionary containing the following:

y0: Sequence of intermediate perturbed states at every iteration. The last element, y0[-1], is then the final perturbed state.
status: The value 0 indicates success, meaning that y0[-1] is an eligible state in the target's basin of attraction. The value 1 indicates that the maximum iteration limit was exceeded.
n_iter: Total number of iterations taken before completion.
time: Total run time of the method (in seconds).
t_int, t_var, t_opt: Lists of run times (in seconds) taken at each iteration for three substeps (integrating the system equations, integrating the variational equation, and forecasting the optimal incremental perturbation, respectively).

**How to invoke *neco*.** Place the file named neco.py in a directory contained in your PYTHONPATH environment variable. From that point, you will be able to import the control routine from the module of the same name:

from neco import neco

…within any python program. Once the differential equations defining the desired network have been implemented as a Python function f (as described in the Procedure section above), the control method can then be invoked within your script according to

output = neco(y0, yt, f, **kwargs)

…where y0 and yt are the initial and desired target states of the network, respectively, and **kwargs represents any (optional) additional arguments that one wishes to specify or change from their default values (constraints, the Jacobian matrix, iteration limits, and so on as described in the Procedure section above). For example, to invoke *neco* with minimum and maximum incremental perturbation sizes of $10^{-3}$ and $10^{-2}$, respectively, and an integration time window of 10 time units, one would enter

output = neco(y0, yt, f, eps0=0.001, eps1=0.01, t_max=10)

...and so on.

As part of this protocol we include the source code (neco_example.py) for two examples of finding control perturbations subject to constraints in a two-dimensional system defined by a particle under the influence of a 1D potential plus dissipation. These control problems are illustrated in Fig. 2 of ref. [1]. These examples are provided mainly for the source code itself (in particular, to give guidance to those unfamiliar with Python on how to define a system of ODEs and call *neco* with various parameters). To actually run the examples, simply invoke python on the example file from the command line

python neco_example.py

**Timing:**

The running time of the algorithm presented here depends on the size of the network. Specifically, in ref. [1], it is demonstrated that the asymptotic scaling of the running time is $O(n^{2.5})$, where n is the number of dynamical variables in the system.

**Anticipated results:**

In ref. [1], we demonstrated our algorithm with a number of example applications, based on the source code implementation included here. These applications involved networks from various domains in the physical and biological sciences. In all cases, we had remarkable success in identifying eligible compensatory perturbations that could control the network under constraints that were appropriate to the particular system in question. In particular, using benchmark networks that were randomly-generated in such a way that compensatory perturbations were guaranteed to exist under the given constraints, our method succeeds in identifying them in 100% of cases.

We note as a final remark that, for clarity, this protocol describes and provides the code only for a simple (albeit, broadly applicable) implementation of our algorithm. It is, however, straightforward to extend and tailor various aspects of the method to suit different but related problems. For example, to

- optimize a different function of the system's trajectory at every iteration (other than distance to the target);
- find an eligible point further in the interior of the desired basin (as opposed to just inside the basin boundary). Such points have the advantage of making this control strategy more robust against noise in the dynamics and uncertainty in the system parameters (we refer to the Supplementary Information of ref. [1]);
- implement subsequent control interventions after monitoring how the system responds to the first one, in a "closed loop" control setting.

This protocol is thus best regarded not as a single algorithm or piece of software, but rather a class of iterative strategies for identifying network control interventions. The choice of the exact goal that these interventions are meant to accomplish in a particular problem is then left up to the experimenter. Details of how to modify the presented approach for the above extensions, as well as others, are beyond the scope of this protocol, but we encourage the reader to contact the authors (at cornelius@u.northwestern.edu or motter@northwestern.edu) for guidance in this area.

**Where to download the source code:**

The source code for NECO and the included example can be downloaded as a combined ZIP file at the following link:

http://www.nature.com/protocolexchange/system/uploads/2647/original/neco_source.zip

**Acknowledgments:**

The authors thank Takashi Nishikawa for illuminating discussions. This work was supported by NSF under Grant DMS-1057128, NCI under Grant 1U54CA143869-01, and a Northwestern-Argonne Early Career Investigator Award to A.E.M.

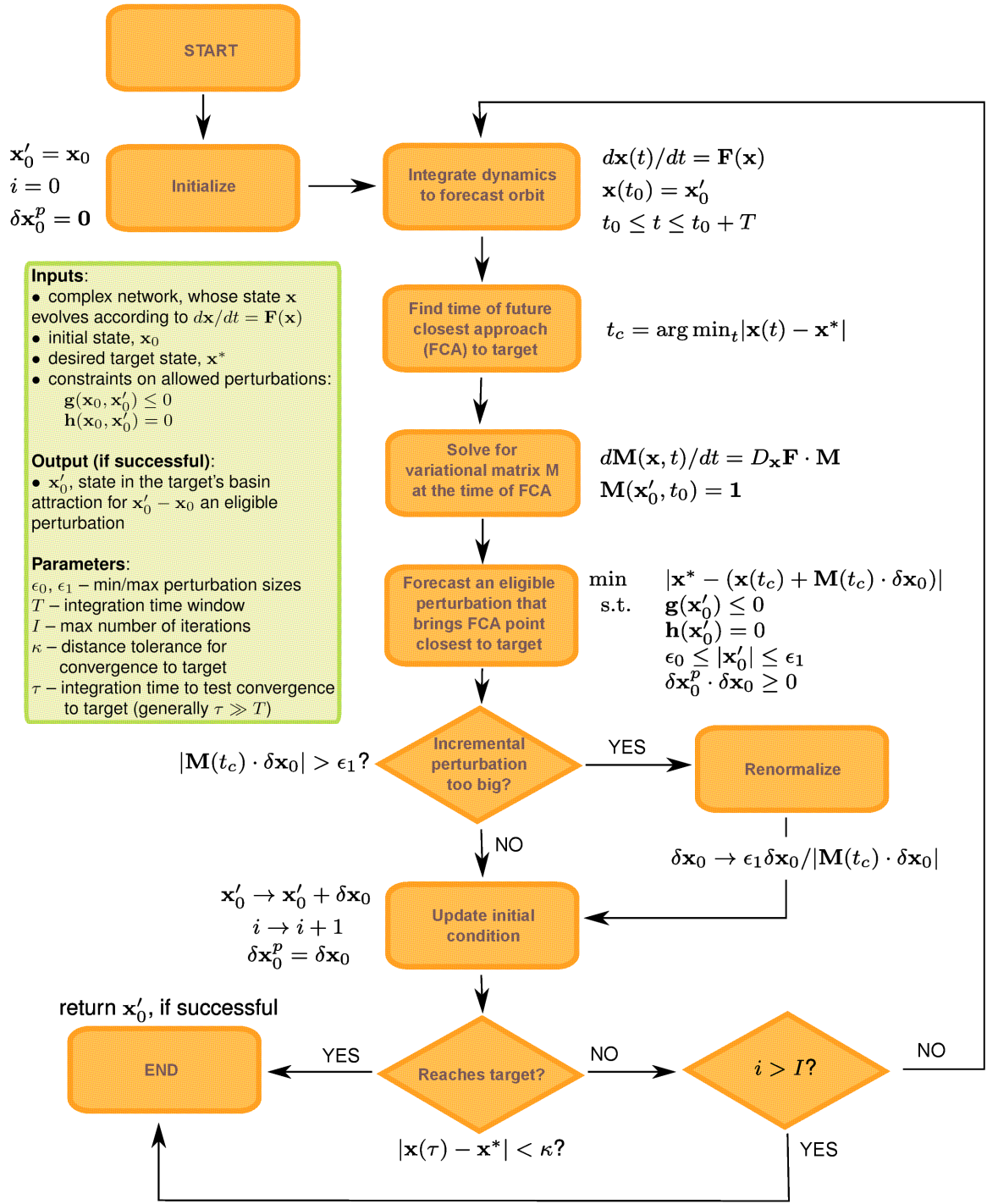

**Fig. 1 Algorithm for identifying eligible control interventions.** The flowchart depicts the steps of the network control algorithm implemented in this protocol. The box on the left describes the inputs, output, and parameters of the algorithm. For the motivation behind each of these parameters, as well as guidance in how to choose their values based on the time and length scales of a particular network, we refer the reader to ref. [1] (the mathematical notation used here is consistent with that used in the reference). For the correspondence between the mathematical objects here and the variables in the software implementation given in this protocol, we refer to the Procedure section.